\documentclass[12pt]{amsart}
\title[Local rigidity of Schottky maps]{
Local rigidity of Schottky maps
}
\author{Sergei Merenkov}
\address{Department of Mathematics\\
University of Illinois\\ 1409 W Green St\\ Urbana, IL
61801\\USA} \email{merenkov@illinois.edu}
\thanks{Supported by NSF grant DMS-1001144}

\date{\today}

\usepackage[active]{srcltx}

\usepackage{amsmath, amssymb, amsthm,latexsym}
\usepackage{graphicx}
\newcommand\C{{\mathbb C}}

\newcommand\N{{\mathbb N}}

\newcommand\Z{{\mathbb Z}}

\newcommand\dee{\partial}

\newcommand\id{\operatorname{id}}

\renewcommand\:{\colon}

\newcommand\no{\noindent}

\newtheorem{theorem}{Theorem}[section]

\newtheorem{corollary}[theorem]{Corollary}

\newtheorem{remark}{Remark}
\newtheorem{lemma}[theorem]{Lemma}

\theoremstyle{definition}

\begin{document}


\abstract
{
We introduce Schottky maps---conformal maps between relative Schottky sets, and study their local rigidity properties. This continues the investigations of relative Schottky sets initiated in~\cite{sM10}. Besides being of independent interest, the latter and current works provide key ingredients in the forthcoming proof of quasisymmetric rigidity of Sierpi\'nski carpet Julia sets of rational functions. 
}
\endabstract

\maketitle

\section{Introduction}\label{s:Intro}
%
%

\no
Let $D$ be a domain (i.e., an open and connected set) in the complex plane $\C$. 
A \emph{relative Schottky set}  $S$ in $D$ is a set 
of the form
$$
S=D\setminus\cup_{i\in I} B_i,
$$ 
where each $B_i$  is an open geometric (round) disc with the closure $\overline B_i$ contained in $D$, and $\overline B_i\cap\overline B_j=\emptyset$ for $i\neq j$. 
Such sets were introduced in~\cite[Section~8]{BKM07}.
The boundary $\dee B_i$ of a disc $B_i$ for each $i\in I$ is referred to as a \emph{peripheral circle}. 
If $D=\C$, a relative Schottky set $S$ is called a \emph{Schottky set}. In this case we allow peripheral circles to degenerate to straight lines, i.e.,  we allow $B_i$ to be an open half-plane for some $i\in I$.

Assume that $S$ is a subset of $\C$ that has an accumulation point $p\in S$, let $U$ be a neighborhood of $p$ in $\C$, and suppose that $f\: S\cap U\to\C$ is a continuous map. We say that $f$ is \emph{conformal} at $p$ if 
$$
\lim_{q\in S,\, q\to p}\frac{f(q)-f(p)}{q-p}
$$   
exists and is non-zero. 
The limit above is called the \emph{derivative} of $f$ at $p$ and it is denoted by $f'(p)$.

Suppose that $D$ and $\tilde D$ are domains contained  in $\C$, and let $S=D\setminus \cup_{i\in I} B_i$ and $\tilde S=\tilde D\setminus \cup_{j\in J}\tilde B_j$ be relative Schottky sets in $D$ and  $\tilde D$, respectively.
Let $U$ be an open subset of $D$ and $f\: S\cap U\to \tilde S$ be a local homeomorphism. We recall that the latter means that for every point $p\in S\cap U$ there exists an open set $V$ that contains $p$ such that $f(S\cap U\cap V)$ is open in $\tilde S$ and the restriction $f|_{V}\: S\cap U\cap V\to f(S\cap U\cap V)$ is a homeomorphism. 
It follows that 
$p\in S\cap U\cap \dee B_i$ for some $i\in I$, if and only if
$f(p)\in \dee \tilde B_j$ for some $j\in J$; see Lemma~\ref{L:Lochomeo} below.
Such a map $f\: S\cap U\to\tilde S$ is called a \emph{Schottky map} if 
it is  conformal at every point $p\in S\cap U$ and
the derivative $f'$ is a continuous function on $S\cap U$.

The existence and, indeed, abundance of Schottky maps is demonstrated in~\cite[Theorems~1.6,~1.2]{sM10}. 
Also, it is a simple observation that if $S$ has measure zero, a Schottky map $f\: S\to \tilde S$ cannot be the restriction of a conformal map from $D$ to $\tilde D$, unless $f$ is the restriction to $D$ of a M\"obius transformation; see~\cite[discussion right after Corollary~1.5]{sM10}.


It is interesting to see how much of the classical theory of analytic functions carries over to Schottky maps. Even the most basic questions are not understood. For example, it is unknown whether every Schottky map between relative Schottky sets of measure zero has higher derivatives in the above sense; see~\cite[Conjecture~1.3]{sM10}. In this note we develop basic properties of Schottky maps and provide, in particular, an analogue of the well known fact that if two analytic functions $f$ and $g$ in a domain $D$ coincide on a set that  has an accumulation point in $D$, then $f=g$ identically in $D$; see Corollary~\ref{C:Un} below. Furthermore, in Section~\ref{S:Appl} we give applications of our results that serve as important steps in a forthcoming joint project with M.~Bonk and M.~Lyubich.  The latter is concerned with the quasisymmetric rigidity of Sierpi\'nski carpet Julia sets of post-critically finite rational functions.

\medskip\noindent 
\textbf{Acknowledgments.} The author thanks  Mario Bonk  for many useful comments that helped to simplify some of the arguments and strengthen the results.

\section{Relative Schottky sets and weak tangents}\label{S:RelSchsets}
\noindent
Below for $p$ and $q\in\C$ we denote by $|p-q|$ the Euclidean distance between these points. 
We denote by $B(p,r)$
an open disc in $\C$ of radius $r>0$ that is centered at $p\in\C$, and  by $\overline B(p,r)$ its closure.

\begin{lemma}\label{L:PC}
Let $S$ be a relative Schottky set in $D$ and $p\in S$. Assume that $r>0$ is such that $B(p,2r)\subseteq D$.
 If $q_0$ and $q_1\in S\cap B(p,r)$, then there exists a rectifiable curve $\gamma$ in $S\cap B(p,2r)$ connecting $q_0$ and $q_1$, and  so that  
$$
{\rm length}(\gamma)\leq \pi |q_1-q_0|.
$$
\end{lemma}
\no
\emph{Proof.}
This is essentially~\cite[Corollary~2.3]{sM10}: we apply~\cite[Lemma~2.1]{sM10} 
to the straight  line $L$ connecting $q_0$ and $q_1$, and use~\cite[Lemma~2.2]{sM10}.
\qed

\medskip

\begin{corollary}\label{C:Lip}
If $S$ is a relative Schottky set in a domain $D$ and $f\: S\cap U\to \tilde S$ is a Schottky map to a relative Schottky set $\tilde S$, then 
$f$ is locally Lipschitz, i.e.,
for each $p\in S\cap U$, there exist  $L\geq0$ and a neighborhood $V$ of $p$, so that 
$$
|f(q)-f(q')|\leq L|q-q'|
$$
for all $q$ and $q'$ in $S\cap V$.
\end{corollary}
\no
\emph{Proof.}
For $p\in S\cap U$ we choose $r>0$  such that $B(p,2r)$ is compactly contained in $D$. If $q$ and $q'$ are in $B(p,r)$, then by Lemma~\ref{L:PC} there exists a rectifiable curve $\gamma$ in $S\cap B(p,2r)$ connecting $q$ and $q'$ with ${\rm length}(\gamma)\leq \pi|q-q'|$. Since the closure $\overline B(p,2r)$ is a compact subset of $D$ and the derivative $f'$ is assumed to be continuous, this derivative is bounded on $S\cap \overline B(p,2r)$, say by $L'\geq0$. 
This gives  
$$
{\rm length}(f(\gamma))\leq L'\, {\rm length}(\gamma).
$$
Because $|f(q)-f(q')|\leq {\rm length}(f(\gamma))$ and ${\rm length}(\gamma)\leq \pi|q-q'|$, we conclude that $f$ is Lipschitz with $L=\pi L'$.
\qed

\medskip

Suppose that  $D$ is a domain in $\C$ and let $S$ be a relative Schottky set in $D$. 
Let $(p_k)$ be a sequence of points in $S$ that converges to $p\in S$. 
A closed set $S_\infty\subseteq \mathbb C$ is called a \emph{weak tangent} of $S$ along the sequence $(p_k)$, if there exists a sequence $(\lambda_k)$ of non-zero complex numbers converging to 0, such that
$$
S_k=(S-p_k)/\lambda_k\to S_\infty\quad {\rm when}\quad k\to \infty,
$$
in the Gromov--Hausdorff sense; see~\cite[Section~8.2]{DS97} for the convergence of subsets of Euclidean spaces. 

A simple application of~\cite[Lemma~8.2]{DS97} gives the following sub-convergence statement.
If $(p_k)$ is a sequence of elements in $S$ such that $\lim_{k\to\infty}p_k=p\in S$, and $(\lambda_k)$ is an arbitrary sequence of non-zero complex numbers with $\lim_{k\to\infty}\lambda_k=0$, then there exists a subsequence $(k_j)$ such that $(S_{k_j})$
converges to $S_\infty$ along $(p_{k_j})$. 
%

A relative Schottky set $S=D\setminus\cup_{i\in I} B_i$ is called \emph{locally porous at} $p\in S$ if there exist a neighborhood $U$ of $p$, a positive   constant $r_0>0$ and a constant $C\geq 1$, such that for every $q\in S\cap U$ and each $r$ with $0<r\leq r_0$, there exists $i\in I$ with $B(q,r)\cap B_i\neq \emptyset$ and
$$
r/C\leq {\rm diam}(B_i)\leq Cr,
$$
where ${\rm diam}(X)$ denotes the diameter of a set $X\subseteq\C$. A relative Schottky set $S$ is called \emph{locally porous} if it is locally porous at every $p\in S$. Every locally porous relative Schottky set has measure zero since it cannot have Lebesgue density points.

\begin{lemma}\label{L:WT}
Assume that $S=D\setminus\cup_{i\in I}B_i$ is a  relative Schottky set that is locally porous at $p\in S$, and let $S_\infty$ be a weak tangent of $S$ along some sequence $(p_k)$ that converges to $p$. Then $S_\infty$ is a locally porous Schottky set in $\C$. Moreover, there exists a constant $C\geq 1$ such that for every $R>0$ there exists a peripheral circle $\dee B$ of $S_\infty$ with $\partial B\cap \overline B(0,R)\neq\emptyset$ and 
$$
R/C\leq {\rm diam}(\dee B)\leq CR.
$$  
\end{lemma}
\no
\emph{Proof.}
The fact that $S_\infty$ is a locally porous Schottky set is essentially proved in~\cite[Lemmas~8.3, 8.4]{BKM07}. The only difference is that in the definition of a weak tangent contained in~\cite{BKM07} it is required that the sequence $(p_k)$ is constant, i.e., that $p_k=p$ for all $k$. 

%
Assume that $S_\infty$ is defined using a sequence $(\lambda_k)$ with  $\lim_{k\to\infty}\lambda_k=0$.
Let $U$ be the neighborhood of $p$, the constants $r_0>0$ and  $C\geq 1$ be those from the definition of local porosity of $S$ at $p$, and let $R>0$ be arbitrary. Then for $k$ large enough, we have $p_k\in S\cap U$ and $0<R|\lambda_k|\leq r_0$. Thus there exists $i\in I$ such that
$B(p_k,R|\lambda_k|)\cap B_i\neq \emptyset$ and 
$$
R|\lambda_k|/C\leq {\rm diam}(B_i)\leq CR|\lambda_k|.
$$
The disc $B_i^k=(B_i-p_k)/\lambda_k$ is a complementary component  of 
$$
S_k=(S-p_k)/\lambda_k.
$$
Also, for every $k$, the disc $B_i^k$ intersects $B(0,R)$ and its diameter satisfies
$$
R/C\leq{\rm diam}(B_i^k)\leq CR.
$$
Since $(S_k)$ converges to $S_\infty$, the sequence $(\dee B_i^k)_k$ converges to a peripheral circle $\dee B$ of $S_\infty$. 
As a limit, the peripheral circle $\dee B$ satisfies
$$
R/C\leq{\rm diam}(\dee B)\leq CR,
$$
and it intersects the closure of the disc $B(0,R)$. This concludes the proof of the lemma.
%
\qed

\medskip

Any M\"obius transformation $f$ that fixes $\infty$ has either the form $f(z)=az+b$ or the form $f(z)=a\bar z+b$, with $a,b\in\C$ and $a\neq 0$. The former holds when $f$ is orientation preserving and the latter when it is orientation reversing.
The following lemma describes the group (under composition) of all M\"obius transformations that fix $\infty$ and preserve a Schottky set $S$, provided the set $S$ contains arbitrarily large peripheral circles. We show that such a group is isomorphic to a discrete subgroup of the semidirect product $\C^\times\rtimes(\Z/2\Z)$, where $\C^\times$ is the multiplicative group of non-zero complex numbers and $\Z/2\Z=\{1,-1\}$ is the group of two elements acting on $\C^\times$ by conjugation. 

%

\begin{lemma}\label{L:Gp}
Suppose that $S$ is a  Schottky set, i.e., a relative Schottky set in $\C$. We assume that it has arbitrarily large peripheral circles, i.e., there exists a sequence of distinct peripheral circles $(\dee B_i)$ of $S$ such that 
$$
\lim_{i\to\infty} {\rm diam}(\dee B_i)=\infty.
$$
Then 
there is no translation that keeps the Schottky set $S$ invariant. 
Moreover, if  the group $G$ of all M\"obius transformations that fix $\infty$ and preserve $S$ contains a non-trivial orientation preserving element, then all the elements of  $G$ share a fixed point in $\C$, and  $G$ is isomorphic to a discrete subgroup of $\C^\times\rtimes(\Z/2\Z)$. Otherwise $G$ is isomorphic to a subgroup of $\Z/2\Z$.
\end{lemma}
\noindent
\emph{Proof.}
Suppose that $f$ is an element of $G$.
It cannot be a translation since, according to our assumption, the Schottky set $S$ contains peripheral circles of arbitrarily large radii. A translation $z\mapsto z+c$ would have to move these large circles, but it clearly 
cannot move the corresponding interior discs off themselves.

We first assume that there exists a non-trivial orientation preserving element $f$ of $G$. Then $f(z)=az+b$ with $a\neq1$, and hence $f$ has a fixed point. Without loss of generality we may assume that the fixed point is 0, and so $f(z)=az$ with $a\neq 1$. 

If $g$ is any other orientation preserving element of $G$, then it has the form $g(z)=c z+d$, where $c\neq 0$. The commutator $[g,f]$ of $g$ and $f$ is
$$
[g,f](z)=g^{-1}f^{-1}gf(z)=z+\frac{d}{c}\left(\frac1{a}-1\right).
$$ 
Since $S$ does not support translations and $a\neq 1$, we must have $d=0$, i.e., we get $g(z)=c z$, and so 0 is also a fixed point of $g$. 

If $g\in G$ is orientation reversing, then $g(z)=c\bar z+d$ with $c\neq 0$. The commutator $[g,f]$ is given by
$$
[g,f](z)=\frac{a}{\bar a}z+\frac{\bar d}{\bar c}\left(\frac1{\bar a}-1\right),
$$
and it is an orientation preserving element of $G$. Since $\bar a\neq 1$, the above implies that we must have $d=0$, and thus $g(z)=c\bar z$. Therefore 0 is again a fixed point of $g$.

We conclude that if $G$ contains a non-trivial orientation preserving element, then all the elements of $G$ share a fixed point in $\C$ and (up to an appropriate change of variables) each $g\in G$ can be written as $g(z)=cz$ or $g(z)=c\bar z$, depending on whether $g$ preserves or reverses the orientation, respectively.
The group $G$ must be discrete because elements of such a group have to preserve the set of peripheral circles. A monomorphism of $G$ into $\C^\times\rtimes(\Z/2\Z)$ is given by $g\mapsto (c,1)$ if $g(z)=cz$, and $g\mapsto(c,-1)$ if $g(z)=c\bar z$. 
%

If the group $G$ is non-trivial but does not have a non-trivial orientation preserving element, then $G$ must be isomorphic to $\Z/2\Z$. Indeed, suppose that  $g\in G$ reverses the orientation. Then the composition $g\circ g$ of $g$ with itself must preserve the orientation, and therefore $g\circ g=\id$. This means that $g^{-1}=g$. Moreover, if $g$ and $h$ are any two orientation reversing elements of $G$, then $g\circ h$ is orientation preserving, and hence $g\circ h=\id$. Thus $h=g^{-1}=g$. 
The group $G$ is therefore isomorphic to $\Z/2\Z$.
\qed

\medskip

\begin{remark}
In this lemma one can replace a Schottky set by a circle domain, i.e., a domain obtained from $\C$ by removing disjoint closed discs. One also has to assume that a circle domain needs to have  arbitrarily large complementary components for the conclusions as in the previous lemma to hold.
\end{remark}


\section{Calculus for Schottky maps}\label{S:Calc}

\no
In this section we prove auxiliary results that will be used in the proof of local uniqueness properties. The main result is Lemma~\ref{L:FTC}, which is a version of the Fundamental Theorem of Calculus for Schottky maps. 

We first establish the following property of local homeomorphisms stated in the Introduction.
\begin{lemma}\label{L:Lochomeo} Assume that $S$ and $\tilde S$ are relative Schottky sets and 
$f\: S\cap U\to \tilde S$ is a local homeomorphism, where $U$ is an open set. Then  
$p\in S\cap U\cap \dee B_i$ for  $i\in I$, if and only if
$f(p)\in \dee \tilde B_j$ for some $j\in J$. Here $\{\dee B_i\}_{i\in I}$ and $\{\dee \tilde B_j\}_{j\in J}$ are the families of peripheral circles of $S$ and $\tilde S$, respectively.
\end{lemma}
\no
\emph{Proof.}
If $p\in S\cap U$, then there exists an open neighborhood $V\subseteq U$ of $p$ such that $f(S\cap V)$ is open in $\tilde S$ and the restriction $f|(S\cap V)\: S\cap V\to f(S\cap V)$ is a homeomorphism. We may assume that $V$ is a small geometric disc centered at $p$ and $f$ restricted to the closure of $S\cap V$ is a homeomorphism onto its image. 

Now, if $p\in \dee B_i$ for $i\in I$, by possibly changing the radius of $V$ we may further assume that the boundary $\dee V$ intersects $\dee B_i$ at exactly two points, and for every other peripheral circle $\dee B_k$ of $S$, the intersection $\dee V\cap \dee B_k$ is either empty or also consists of exactly two points.  The set $S\cap\dee V$ combined with the arcs $\dee B_k\cap V$ of all the peripheral circles $\dee B_k$ of $S$ (including $\dee B_i$) that intersect $\dee V$ form a Jordan curve $\alpha$ contained in the closure of $S\cap V$.  
This follows from the assumption that peripheral circles are disjoint geometric circles.
The image $f(\alpha)$ is thus a Jordan curve in $\tilde S$ that contains $f(p)$. 

If $f(p)$ does not belong to any of the peripheral circles of $\tilde S$, then the Jordan--Schoenflies theorem implis that $f(p)\in f(\alpha)$ is an accumulation point for the set $\tilde S\setminus f(S\cap V)$. Once proved, this immediately leads to a contradiction because the set $\tilde S\setminus f(S\cap V)$  is closed in $\tilde S$ and it does not contain $f(p)$. 

To prove the above statement, one first observes that according to the Jordan--Schoenflies theorem any open disc centered at $f(p)$ must contain points in the connected component $C$ of the complement of $f(\alpha)$ that does not contain points of $f(S\cap V)$. 
Let $(r_k)$ be a sequence of positive numbers converging to 0, and for each $k$ let $D_k$ denote the open disc centered at $f(p)$ whose radius is $r_k$. 
If $p_k$ is a point in $D_k\cap C$, then starting from some number $K_0\in\N$, each $p_k,\ k\ge K_0$,  belongs to the domain $\tilde D$ that contains $\tilde S$. If there exists a subsequence $(r_{k_l})$ such that each $p_{k_l}$ is in $\tilde S$, then we are done. Otherwise, for some $K_1\in\N$, each $p_k,\ k\geq K_1$, belongs to the disc bounded by a peripheral circle $\dee \tilde B_{j_k}$ of $\tilde S$. If all $\dee \tilde B_{j_k}$ are the same peripheral circle for $k\ge K_2$, for some $K_2\in\N$, then $f(p)$ belongs to this peripheral circle, a contradiction.  Hence there exists a subsequence, which we assume to be the whole sequence $(\dee \tilde B_{j_k})_{k\ge K_1}$, such that $\dee \tilde B_{j_{k'}}\neq \dee \tilde B_{j_{k''}}$ for $k'\neq k''$. For each $k\ge K_1$, the intersection $\dee \tilde B_{j_k}\cap C$ is non-empty because the boundary $f(\alpha)$ of $C$ is a Jordan curve. By choosing $p_k'\in\dee \tilde B_{j_k}\cap C$, we obtain a new sequence $(p_k')$ that converges to $f(p)$ and such that $p_k'\in \tilde S\cap C\subseteq \tilde S\setminus f(S\cap V)$. This means that $f(p)$ is an accumulation point for the set $\tilde S\setminus f(S\cap V)$ and the proof of one implication is complete.

The converse statement, namely that if $f(p)$ belongs to $\partial \tilde B_j$ for $j\in J$, implies $p\in \partial B_i$ for some $i\in I$, is proved using the same argument applied to the local inverse of $f$. 
\qed

\medskip

In the following, when we say that $f\: S\cap U\to\tilde S$ is a Schottky map, then it is understood that $S$ and $\tilde S$ are relative Schottky sets in domains $D$ and $\tilde D$, respectively, that $U$ is an open subset of $D$, and that $f$ is conformal at every point of $S\cap U$ with $f'$ being a continuous function on $S\cap U$.

\begin{lemma}\label{L:Comp}
Let $f\: S\cap U\to S'$ and $g\: S'\cap V\to \tilde S$ be Schottky maps such that $f(S\cap U)\subseteq S'\cap V$. Then $h=g\circ f\: S\cap U\to \tilde S$ is a Schottky map and at each $p\in S\cap U$ we have
$$
h'(p)=g'(f(p))f'(p).
$$
\end{lemma}
\no
\emph{Proof.}
This follows from the linear approximations of the maps under consideration.
\qed

\medskip
\begin{lemma}\label{L:Inv}
Let $f\: S\cap U\to\tilde S$ be a Schottky map. Suppose that $p\in S\cap U$ is an arbitrary point  and $g$ is the local inverse of $f$ defined in the intersection of a neighborhood of $f(p)$ with $\tilde S$. Then $g$ is conformal at $f(p)$ and 
$$
g'(f(p))=1/f'(p).
$$ 
\end{lemma}
\no
\emph{Proof.}
This follows immediately  from the definition of derivative. 
\qed

\medskip

\begin{corollary}\label{C:biLip}
If $f\:S\cap U\to\tilde S$ is a Schottky map, it is locally bi-Lipschitz, i.e.,  
for each $p\in S\cap U$, there exist  $L\geq 1$ and neighborhoods $V$ and $\tilde V$ of $p$ and $f(p)$, respectively, so that $f\: S\cap V\to \tilde S\cap \tilde V$ is a homeomorphism satisfying
$$
\frac1{L}|q-q'|\leq|f(q)-f(q')|\leq L|q-q'|
$$
for all $q$ and $q'$ in $S\cap V$.
\end{corollary}
\no
\emph{Proof.}
From Corollary~\ref{C:Lip} we know that $f$ is locally Lipschitz. Also, it follows from Lemma~\ref{L:Inv} that the local inverse $g$ of the local homeomorphism $f$ is a Schottky map defined in a neighborhood of $f(p)$. Hence by the same Corollary~\ref{C:Lip} the map $g$ is Lipschitz in some neighborhood of $f(p)$. We therefore conclude that $f$ is locally bi-Lipschitz.
\qed

\medskip

\begin{lemma}\label{L:FTC}
Let 
$f\:S\cap U\to\tilde S$ be a Schottky map. If $q_0, q_1\in S\cap U$, and $\gamma$ is a rectifiable curve in $S\cap U$ with initial point $q_0$ and terminal point $q_1$, 
then
$$
f(q_1)-f(q_0)=\int_\gamma f'(z)dz.
$$
\end{lemma}
\no
\emph{Proof.}
Let $l$ be the length of $\gamma$ and $\gamma\:[0, l]\to S\cap U$ be an arc-length parametrization. 
We set $h(t)=f(\gamma(t))$. Since $\gamma$ is parametrized by arc length and $f$ is locally Lipschitz on a compact set $\gamma$, the map $h\:[0,l]\to\tilde S$ is absolutely continuous. Moreover, $h'(t)=f'(\gamma(t))\gamma'(t)$ for almost  every $t\in[0,l]$. 
Therefore,
$$
f(q_1)-f(q_0)=h(l)-h(0)=\int_0^l h'(t)dt
=\int_\gamma f'(z)dz.
\qed
$$

\medskip

%

The following will be used repeatedly in the rest of the paper.

Given $L\ge1$, let $(f_k)$ be a sequence where each $f_k\: S\cap U\to\tilde S$ is a Schottky map that is $L$-bi-Lipschitz onto its image.  Suppose that $(p_k)$ is a sequence of elements in $S\cap U$ with $\lim_{k\to\infty}p_k=p\in S\cap U$, and 
let $(\lambda_k)$ be a sequence of non-zero complex numbers with $\lim_{k\to\infty}\lambda_k=0$. Further,
assume that $(\tilde p_k)$ is a sequence of elements in $\tilde S$ that converges to $\tilde p\in\tilde S$. 
Finally, define
$$
F_k(z)=\frac{f_k(p_k+\lambda_k z)-\tilde p_k}{\lambda_k},
$$
a map from $S_k(U)=(S\cap U-p_k)/\lambda_k$ to $\tilde S_k=(\tilde S-\tilde p_k)/\lambda_k$.
By possibly passing to  subsequences, we may assume that 
$(S_k(U))$ converges to a weak tangent $S_\infty$ and $(\tilde S_k)$ converges to a weak tangent $\tilde S_\infty$.
 The maps $F_k$  are $L$-bi-Lipschitz on compacta for all  $k$ large enough. This means that, given any compact set $K\subseteq\C$, there exists $N\in \N$, so that for each $k\ge N$ the map $F_k$ restricted to $S_k(U)\cap K$ is $L$-bi-Lipschitz onto its image.
Hence, a certain subsequence $(F_{k_j})$ of $(F_k)$ converges  to a bi-Lipschitz map $F_\infty\colon S_\infty\to \tilde S_\infty$; see~\cite[Lemmas~8.6 and 8.7]{DS97}. Note that the maps $F_k$ are not defined on the same set, but for each $k$, the map $F_k$ is defined on $S_k(U)$, and the sequence $(S_k(U))$ converges to $S_\infty$. The  convergence of the sequence $(F_{k_j})$ here means that for every $q\in S_\infty$ and every sequence $(q_j)$ with $q_j\in S_{k_j}(U)$ and $\lim_{j\to\infty}q_j=q$, we have $\lim_{j\to\infty} F_{k_j}(q_j)=F_\infty(q)$; see~\cite[Definition~8.4]{DS97}.

If we assume, in addition, that the set $S$ is locally porous at $p$, then the set $S_\infty$ is locally porous according to Lemma~\ref{L:WT}. In particular, it has measure zero. Therefore, by~\cite[Theorem~1.2]{BKM07}, the map $F_\infty$ is the restriction to $S_\infty$ of a M\"obius transformation. Clearly $F_\infty$ must fix $\infty$. 

An important special case of the previous discussion is when each $f_k$ is the same Schottky map $f$. By Corollary~\ref{C:biLip}, the map $f$ is bi-Lipschitz when restricted to some neighborhood of each $p\in S\cap U$.


\section{Local uniqueness of Schottky maps}\label{S:Locuniq}
\no
This section contains local rigidity results for Schottky maps.

\begin{theorem}\label{T:Un}
Assume that $S$ is a locally porous relative Schottky set in a domain $D\subseteq\C$, and let $U\subseteq D$ be an open  set such that $S\cap U$ is connected. 
Let $f\: S\cap U\to S$ be a Schottky map and $p\in S\cap U$. Assume further that $f(p)=p$ and $f'(p)=1$. Then  $f=\id$ in $S\cap U$. 
\end{theorem}
\no
\emph{Proof.}
Suppose  first that there exists a sequence $(p_k),\ p_k\in S\cap U$, that converges to $p$ and such that $f(p_k)\neq p_k$. 
Let $\lambda_k=f(p_k)-p_k$.  Each $\lambda_k\neq0$, and  $\lim_{k\to\infty}\lambda_k=0$ because $f$ is continuous and $f(p)=p$. Let us consider 
$$
F_k(z)=\frac{f(p_k+\lambda_k z)-p_k}{\lambda_k}=\frac{f(p_k+\lambda_k z)-f(p_k)}{\lambda_k}+1,
$$ 
a map from $S_k(U)=(S\cap U-p_k)/\lambda_k$ to $S_k=(S-p_k)/\lambda_k$.

There exists a subsequence $(k_j)$ such that the sequences
$(S_{k_j}(U))$ and $(S_{k_j})$  converge to a weak tangent space $S_\infty$. Without loss of generality we assume that $(k_j)$ is the whole sequence $(k)$.
Let $q$ be a point in $S_\infty$. Then for each $k$ there exists $q_k\in S_k(U)$ with $\lim_{k\to\infty}q_k=q$.

For large $k$, let $\gamma_k$ be a rectifiable curve in $S\cap U$ that connects $p_k$ and $p_k+\lambda_k q_k$ and has length $\le\pi|q_k\lambda_k|$, as guaranteed by  Lemma~\ref{L:PC}. Since $f\: S\cap U\to S$ is a Schottky map, Lemma~\ref{L:FTC} gives
$$
\frac{f(p_k+\lambda_kq_k)-f(p_k)}{\lambda_k}=\frac1{\lambda_k}\int_{\gamma_k}f'(z)dz.
$$ 
Because $f'(p)=1$ and $f'$ is continuous at $p$, for every $\epsilon>0$ there exists $N\in \N$ such that for every $k\geq N$ we have
$$
\sup\{|f'(z)-1|\: z\in \gamma_k\}<\epsilon.
$$ 
Therefore, since $\int_{\gamma_k}1dz=\lambda_kq_k$, we have
$$
\begin{aligned}
|F_k(q_k)-1-q_k|&=\bigg|\frac{f(p_k+\lambda_kq_k)-f(p_k)}{\lambda_k}-q_k\bigg|\\
&=\bigg|\frac1{\lambda_k}\int_{\gamma_k}\big(f'(z)-1\big)dz\bigg|\\
&\leq\epsilon\, {\rm length}(\gamma_k)/|\lambda_k|
\leq \pi|q_k|\epsilon,\quad k\geq N.
\end{aligned}
$$
Since $\lim_{k\to\infty}q_k=q$, we conclude that $\lim_{k\to\infty}F_k(q_k)=q+1$. 

The sequence $(F_k)$ consists of maps uniformly bi-Lipschitz on compacta for all $k$ large enough. Thus a certain subsequence of $(F_k)$ converges  to a map $F_\infty\colon S_\infty\to S_\infty$.
The above shows that $F_\infty(q)=q+1$ for every $q\in S_\infty$.
On the other hand, by Lemma~\ref{L:WT}, the weak tangent $S_\infty$ contains arbitrarily large peripheral circles, and hence, by Lemma~\ref{L:Gp}, it cannot support translations. 
This is a contradiction.

The above argument implies that the set of all points $p$ in $S\cap U$ such that $f(p)=p$ and $f'(p)=1$ is open in $S\cap U$. Since $f$ and $f'$ are continuous, it is also a closed subset of $S\cap U$. But $S\cap U$ is assumed to be connected. Therefore, $f$ is the identity map. 
\qed

\medskip

\begin{corollary}\label{C:Un}
Let $S$ be a locally porous relative Schottky set in $D\subseteq \C$,  and suppose that $U\subseteq D$ is an open  set such that $S\cap U$ is connected. 
Let $f$ and $g\: S\cap U\to \tilde S$ be Schottky maps into a relative Schottky set $\tilde S$ in a domain $\tilde D$, and consider 
$$
E=\{p\in S\cap U\: f(p)=g(p)\}.
$$  
Then $E=S\cap U$, provided $E$ has an accumulation point in $U$. 
\end{corollary}
\no
\emph{Proof.}
Since $f$ and $g$ are continuous, the set $E$ is a closed subset of $S\cap U$. 

We assume that $E$ has an accumulation point in $U$, and let $p\in S\cap U$ be such a point. We consider the map $h=g^{-1}\circ f$ defined in a neighborhood $V$ of $p$. 
The map $h\: S\cap V\to S$ is a Schottky map as follows from Lemmas ~\ref{L:Comp} and \ref{L:Inv}. Also, since $p$ is an accumulation point for the set $\{q\in S\cap V\: h(q)=q\}$, we have $h'(p)=1$. Theorem~\ref{T:Un} implies that $h=\id$ in the connected component of $S\cap V$ that contains $p$. Hence $f=g$ in some neighborhood of $p$. This shows that the set of accumulation points $E'$ of $E$ is open in $S\cap U$.  Clearly $E'$ is also closed in $S\cap U$. Since $S\cap U$ is connected, we have $E'=S\cap U$. But $E'\subseteq E\subseteq S\cap U$, and thus $E=S\cap U$. 
\qed

\medskip

\begin{corollary}\label{C:UUn}
Let $S$ be a locally porous relative Schottky set in a domain $D\subseteq\C$, and  $U\subseteq D$ be an open set such that $S\cap U$ is connected. 
Assume that $f$ and $g\: S\cap U\to \tilde S$ are Schottky maps into a relative Schottky set $\tilde S$ in $\tilde D$, and let $p\in S\cap U$. If $f(p)=g(p)$ and $f'(p)=g'(p)$, then  $f=g$ in $S\cap U$. 
\end{corollary}
\no
\emph{Proof.}
We consider the map $h=g^{-1}\circ f$ defined in a neighborhood $V$ of $p$. By replacing $S\cap V$ by its connected component containing $p$, we may assume that $S\cap V$ is connected. The map $h\: S\cap V\to S$ is a Schottky map by Lemmas ~\ref{L:Comp} and \ref{L:Inv}, and we have $h(p)=p$ and $h'(p)=1$. Theorem~\ref{T:Un} implies that $h=\id$ in $S\cap V$, and so $f=g$ in a neighborhood of $p$. We can now apply Corollary~\ref{C:Un} to conclude that $f=g$ in $S\cap U$.
\qed

\medskip

Assume that $S$ is a locally porous relative Schottky set in a domain $D\subseteq\C$,
a point  $p\in S$ is arbitrary, and $f\: S\cap U\to S$ is a Schottky map that is defined in some neighborhood $S\cap U$ of $p$,  such that $f(p)=p$. We identify two such maps $f$ and $g$  if there exists a neighborhood $V$ of $p$ so that $f|_{S\cap V}=g|_{S\cap V}$. 
This is an equivalence relation, and let  $G$ denote the group whose elements are equivalence classes $[f]$. The group operation is defined by $[f]\cdot[g]=[f\circ g]$. 

\begin{corollary}\label{C:Sm}
 The group $G$ is isomorphic to a discrete subgroup of $\C^\times$. A monomorphism of $G$ into $\C^\times$  is given by $[f]\mapsto f'(p)$, where 
 $f'(p)\in \C^\times$ is the derivative of $f$ at $p$.
\end{corollary}
\no
\emph{Proof.}
Suppose that $S_\infty$ is a weak tangent space of $S$ at $p$. By Lemma~\ref{L:WT}, the Schottky set $S_\infty$ is locally porous and has arbitrarily large peripheral circles.  Every Schottky map $f\in G$ induces an orientation preserving  M\"obius transformation $\Lambda$ that fixes $\infty$ and leaves $S_\infty$ invariant, and it is given by $\Lambda(z)=f'(p) z$.
According to Lemma~\ref{L:Gp}, the group $G'$ of such maps $\Lambda$ is isomorphic to a discrete subgroup of $\C^\times\rtimes\Z/2\Z$. In fact, it is isomorphic (via the map given by $\Lambda\mapsto f'(p)$) to a discrete subgroup of $\C^\times$ because  the maps in $G'$ are orientation preserving.

We now look at the map $\phi\: G\to G'$ given by
$$
\phi([f])=\Lambda.
$$
Lemma~\ref{L:Comp} tells us that this is a homomorphism, and Corollary~\ref{C:UUn} gives that it is injective. The corollary follows.
\qed

\medskip

Any discrete subgroup of $\C^\times$ is isomorphic to $(\Z/n\Z)\times\Z$, where $n\in\N$ and the factor $\Z$ may be missing. In other words, every element of such a discrete subgroup can be decomposed in a unique way as the product of an element in $\Z/n\Z$ and an element in $\Z$. Elements of $\Z/n\Z$ are rotations about the origin, and those of $\Z$ are scalings by an integer power of some $\lambda\in\C^\times$  with $|\lambda|\neq 1$.

The next lemma along with Corollary~\ref{C:Sm} show that for some $n\in\N$ there is an at most $(\Z/n\Z)\times\Z$ supply of the equivalence classes of Schottky maps that agree at a given point $p$,  where again two Schottky maps are equivalent if they agree in some neighborhood of $p$.

\begin{lemma}\label{L:SmallSupply}
Let $S$ be a locally porous relative Schottky set in $D$ and suppose that $\tilde S$ is a relative Schottky set in $\tilde D$. We assume that $U\subseteq D$ is an open set and let $p\in S\cap U$. Then  if 
 $f$ and  $g\:S\cap U\to\tilde S$ are two Schottky maps with $g(p)=f(p)$, we have  $g=f\circ h$ in $S\cap V$, for some open neighborhood $V$ of $p$ and a Schottky map $h$ with $h(p)=p$, so that $[h]$ is an element of the group $G$ defined above.
\end{lemma}
\no
\emph{Proof.}
We set $h=f^{-1}\circ g$ in a neighborhood of $p$ where this composition is defined and homeomorphic. This is a Schottky map with $h(p)=p$. Thus $[h]$ is an element of the group $G$ and $g=f\circ h$ in $S\cap V$ for some open neighborhood $V$ of $p$. 
\qed

\section{Applications}\label{S:Appl}

\no
In this section we derive consequences of the previous results that play a key role in a forthcoming joint project with M.~Bonk and M.~Lyubich. As mentioned above, the project concerns quasisymmetric rigidity of Sierpi\'nski carpet Julia sets of post-critically finite rational functions. 

The following lemma shows that uniformly bi-Lipschitz Schottky maps near  a locally porous point $p$ either fix this point or move it a definite amount, provided there is a  non-trivial map fixing $p$.

\begin{lemma}\label{L:Fp} 
%
Let $S$ be a relative Schottky set in a domain $D\subseteq \C$, suppose that $p\in S$ and $U$ is an open neighborhood of $p$. We assume that $S$ is locally porous at $p$. 
Further, we assume that there exists a Schottky map 
$f\: S\cap U\to S$ with $f(p)=p$ and $f'(p)\neq1$. Then for every open neighborhood $V$ of $p$ and every $L\geq1$, there exists $\delta>0$ such that if $h\:S\cap V\to S$ is an $L$-bi-Lipschitz Schottky map onto its image, then either $h(p)=p$ or else $|h(p)-p|\geq\delta$.
\end{lemma}
%
%
%
\no
\emph{Proof.}
Assume for contradiction that there exist an open neighborhood $V$ of $p$, a constant $L\geq 1$, and a sequence $(h_k)_{k\in\N}$, where each $h_k\: S\cap V\to S$ is a $L$-bi-Lipschitz Schottky map onto its image, such that $h_k(p)\neq p$ and $\lim_{k\to\infty}h_k(p)=p$. 
Let $\lambda_k=h_k(p)-p$ for $k\in\N$. Then $\lambda_k\neq 0$ for $k\in\N$, and $\lim_{k\to\infty}\lambda_k=0$. We consider
$$
H_k(z)=\frac{h_k(p+\lambda_k z)-p}{\lambda_k},
$$
a map from $S_k(V)=(S\cap V-p)/\lambda_k$ to $S_k=(S-p)/\lambda_k$. 
There exists a subsequence $(k_j)$ such that the sequences $(S_{k_j}(V))$ and $(S_{k_j})$ converge to a weak tangent $S_\infty$.  
By Lemma~\ref{L:WT}, the weak tangent $S_\infty$ 
a locally porous Schottky set that contains arbitrarily large peripheral circles.


By possibly passing to a further subsequence, we may assume that $(H_k)$ converges to a M\"obius map $H\: S_\infty\to S_\infty$ that fixes $\infty$.  
Since $H_k(0)=1$ for all $k\in\N$, we have either $H(z)=az+1$ or $H(z)=a\bar z+1$, with $a\neq 0$.
Also, since $f\: S\cap U\to S$ is conformal and $f(p)=p$, this map induces a map  $\Lambda\: S_\infty\to S_\infty$ given by $\Lambda(z)=\lambda z$, where $\lambda=f'(p)$ and thus $\lambda\neq1$. This contradicts Lemma~\ref{L:Gp} since $H$ and $\Lambda$ have no common fixed point in $\C$.
\qed

\medskip

%


The last theorem below  shows that a convergent sequence of Schottky homeomorphisms from a fixed domain into a fixed target stabilizes, provided  that  there exists a non-trivial Schottky map fixing a point in the domain. 

\begin{theorem}\label{T:Der}
Let $S$ be a locally porous relative Schottky set in a  domain $D\subseteq\C$, and $p\in S$ be an arbitrary point. 
Suppose that $U\subseteq D$ is an open neighborhood of $p$ such that $S\cap U$ is connected. We assume that there exists a Schottky map 
$f\: S\cap U\to S$ with $f(p)=p$ and $f'(p)\neq1$.
Let $\tilde S$ be a relative Schottky set in a domain $\tilde D$ and let $(h_k)_{k\in\N}$ be a sequence of Schottky maps $h_k\: S\cap U\to \tilde S$.  We assume that for each $k\in\N$ there exists an open set $\tilde U_k$ so that the map $h_k\: S\cap U\to \tilde S\cap \tilde U_k$ is a homeomorphism, and the sequence $(h_k)$
converges locally uniformly to a  homeomorphism $h\:S\cap U\to \tilde S\cap\tilde U$, where $\tilde U$ is an open set. 
Then there exists $N\in \N$ such that $h_k=h$ in $S\cap U$ for all $k\geq N$.
\end{theorem}
\no
\emph{Proof.}
%
We first assume that $p$ does not belong to any of the peripheral circles of $S$.
As usual, we write $S=D\setminus\cup_{i\in I}B_i$, where $B_i$ is an open geometric disc with closure contained in $D$, and $\overline B_i\cap \overline B_j=\emptyset$ if $i\neq j$.

Consider a decomposition of the plane whose elements are the closures $\overline B_i$ for all $i\in I$, as well as points of $\C$ that do not belong to any of these closed discs. This is an upper semicontinuous decomposition of the the plane into continua that do not separate the plane; see~\cite{rM25}.
By Moore's theorem~\cite{rM25}, the corresponding decomposition space $\mathcal D$, i.e., the topological space obtained by collapsing each $\overline B_i$ to a point, is homeomorphic to the plane. We therefore identify $\mathcal D$ with the plane. If $P\in\mathcal D$ corresponds to $p$, in any neighborhood of $P$ one can find a circle $C$ centered at $P$ that does not pass through any of the points that correspond to   $\overline B_i,\ i\in I$. This is possible because there are only countably many elements in $I$. If the neighborhood of $P$ is chosen to be small enough,
the Jordan curve $c$ that corresponds to $C$ under the decomposition is contained  in $S\cap U$ and has the following properties. It does not intersect any of the peripheral circles of $S$;  the Jordan domain $\Omega$ bounded by $c$ is contained in $U$ and contains the point $p$. 

Since the sequence $(h_k)$ converges uniformly in $S\cap \overline \Omega$ to a homeomorphism, \cite[Theorem~1.4]{sM10} implies that the maps $h_k$ are uniformly bi-Lipschitz in a neighborhood of $p$. Therefore the  Schottky maps $h_{k+1}^{-1}\circ h_k$ are also uniformly bi-Lipschitz in a (perhaps smaller) neighborhood of $p$. Because $h_{k+1}^{-1}\circ h_k(p)$ converges to $p$, Lemma~\ref{L:Fp} implies that $h_k(p)=h_{k+1}(p)=\tilde p\in\tilde S$ for $k$ large enough.

We now claim that there are only finitely many 
distinct maps $h_k$.
Indeed, 
since 
the maps $h_k$ are uniformly bi-Lipschitz in a neighborhood of $p$, there exists $L\geq 1$ such that
$$
\frac1{L}\leq |h_k'(p)|\leq{L}
$$
for all $k$. By Lemma~\ref{L:SmallSupply}, for $j$ and $k$ large enough,  the maps $h_j$ and $h_{k}$ differ in some neighborhood of $p$ by a composition with a Schottky map $g$, where $[g]$ an element of the group $G$. The latter is isomorphic to a subgroup of $(\Z/n\Z)\times\Z$ for some $n\in\N$. The elements of the factor $\Z$, if present, correspond to scalings on the level of weak tangents. The above bounds for $|h_k'(p)|$ imply that there can only be finitely many distinct values $h_k'(p)$ such that $|h_k'(p)|\neq 1$, and hence the total number of possible values for $h_k'(p)$ is finite.
 Therefore  the assumption that $S\cap U$ is connected and Corollary~\ref{C:UUn} imply that there can only be finitely many distinct maps $h_k$, as claimed.

Finally, since the sequence $(h_k)$ converges to $h$, we conclude that, in the case when $p$ does not belong to any of the peripheral circles of $S$, there exists $N\in\N$ such that $h_k=h$ in $S\cap U$ for all $k\geq N$.

Now suppose that $p$ belongs to a peripheral circle of $S$, say $\dee B_{i_0}$, for some $i_0\in I$. We are going to reduce this case to the previous one so that we can again apply~\cite[Theorem~1.4]{sM10}.
We consider a decomposition of the plane whose elements are $\overline B_i$ for all $i\in I\setminus\{i_0\}$, as well as points that do not belong to any of these discs. This is again an upper semicontinuous decomposition of $\C$ by continua that do not separate the plane. The corresponding decomposition space, still denoted by $\mathcal D$, is thus homeomorphic to the plane. So we identify $\mathcal D$ with the plane, and, moreover, using the Jordan--Schoenflies theorem we may assume that the Jordan curve in $\mathcal D$ that corresponds to $\dee B_{i_0}$ is a geometric circle $A$. 

The point $P$ that corresponds to $p$ belongs to $A$. Similar to the above, in any small neighborhood of $P$, we can find a circle $C$ centered at $P$ that does not contain any  of the points that correspond to $\overline B_i,\ i\in I\setminus\{i_0\}$. Let $C'$ denote the closed arc of the circle $C$ whose end points belong to $A$ and so that $C'$ is contained in the complement of the disc bounded by $A$. 
If the neighborhood of $P$ is chosen small enough,
the Jordan arc $c'$ that corresponds to $C'$ under the decomposition satisfies the following properties. The arc $c'$ is contained  in $S\cap U$ and it does not intersect any of the peripheral circles of $S$ other than $\dee B_{i_0}$. Also, if $b$ denotes the arc of $\dee B_{i_0}$  that contains $p$ and has the same end points as $c'$, the Jordan domain $\Omega$ that is bounded by $c'$ and $b$  satisfies $\Omega\subseteq U$.

We can now use reflection in $b$ to get a relative Schottky set $S'$ in the Jordan domain $\Omega'$ whose boundary is the Jordan curve that consists of the arc $c'$ as well as  the arc $c''$, the reflection of $c'$ in $b$. 
Because the image of $b$ under each $h_k$ is an arc on a peripheral circle of $\tilde S$, 
by using the Schwarz reflection principle we can extend each map $h_k$ into the closure of $S'$.   The extension of each $h_k$ is again a Schottky map. 
Moreover, for $k$ large enough, all the maps $h_k$ take $b$ to arcs on the same peripheral circle $\dee\tilde B_{j_0}$ of $\tilde S$. This holds because the sequence $(h_k)$
converges locally uniformly in $S\cap U$ to a  homeomorphism $h$. 
Thus we are in the setup of the previous case 
where the target relative Schottky set $\tilde S$ should be replaced by its double across $\dee\tilde B_{j_0}$.
The theorem follows.
\qed

\end{document}